\documentclass{amsart}
\usepackage[margin=1in]{geometry}
\usepackage{amsmath,amssymb,amsthm,enumitem,mathtools,hyperref}
\usepackage[numbers]{natbib}
\newtheorem{thm}{Theorem}[section]
\newtheorem{prop}[thm]{Proposition}
\newtheorem{cor}[thm]{Corollary}
\newtheorem{lem}[thm]{Lemma}
\theoremstyle{definition}
\newtheorem{defn}[thm]{Definition}
\newtheorem{example}[thm]{Example}
\newtheorem{remark}[thm]{Remark}
\newcommand{\N}{\mathbb{N}}
\newcommand{\Z}{\mathbb{Z}}
\DeclareMathOperator{\Ap}{Ap}
\newcommand{\ints}[1]{\left[#1\right]}
\DeclareMathOperator{\PP}{\mathrm{P}}
\DeclareMathOperator{\CC}{\mathrm{C}}
\DeclareMathOperator{\nn}{\mathrm{n}}
\DeclareMathOperator{\HH}{\mathrm{H}} % set of gaps
\DeclareMathOperator{\mult}{\mathrm{m}} % multiplicity
\DeclareMathOperator{\genus}{\mathrm{g}} % genus
\DeclareMathOperator{\embdim}{\mathrm{e}} % embedding dimension

\makeatletter
\let\@@pmod\pmod
\DeclareRobustCommand{\pmod}{\@ifstar\@pmods\@@pmod}
\def\@pmods#1{\mkern4mu({\operator@font mod}\mkern 6mu#1)}
\makeatother

\title{Equidistribution of numerical semigroup gaps modulo $m$}
\author{Caleb McKinley Shor}
\address{Department of Mathematics, Western New England University, Springfield, MA 01119}
\email{cshor@wne.edu}
\subjclass[2010]{11D04, 20M14}
\keywords{equidistributed multisets modulo $m$, numerical semigroups, gaps, generalized arithmetic sequences, Ap\'ery sets}
\date{\today}

\begin{document}
\begin{abstract}
For a positive integer $m$, a finite set of integers is said to be equidistributed modulo $m$ if the set contains an equal number of elements in each congruence class modulo $m$. In this paper, we consider the problem of determining when the set of gaps of a numerical semigroup $S$ is equidistributed modulo $m$. Of particular interest is the case when the nonzero elements of an Ap\'ery set of $S$ form an arithmetic sequence. We explicitly describe such numerical semigroups $S$ and determine conditions for which the sets of gaps of these numerical semigroups are equidistributed modulo $m$.
\end{abstract}
\maketitle

\section{Introduction}\label{sec:intro}
Let $\N_0$ denote the set of non-negative integers, a monoid under addition, and let $\N$ denote the set of positive integers. A \emph{numerical semigroup} $S$ is an additive submonoid of $\N_0$ with finite complement. Elements of the complement are called \emph{gaps}. Any numerical semigroup can be written as the set of all finite non-negative linear combinations of a finite set $G=\{g_1,\dots,g_k\}\subset\N$ with $\gcd(G)=1$. We denote this by $S=\langle G\rangle$ or $S=\langle g_1,\dots,g_k\rangle$.

In \cite{WangWang2008}, Wang \& Wang studied alternate Sylvester sums for numerical semigroups of the form $S=\langle a,b\rangle$ with $\gcd(a,b)=1$. Among their results are formulas for the numbers of even gaps and odd gaps of $S$. In particular, they found that there are as many even gaps as there are odd gaps precisely when $a$ and $b$ are both odd. Put another way, they determined when $S$ has an equal number of gaps in each congruence class modulo 2.

In this paper, we consider the following problem: for a general numerical semigroup $S$ and $m\in\N$, does $S$ have an equal number of gaps in each congruence class modulo $m$? When this occurs, we say that the set of gaps of $S$ is \emph{equidistributed modulo $m$}.

Our primary result (Proposition \ref{prop:ED iff congruence}) states that for any nonzero $a\in S$ with $\gcd(a,m)=1$, the set of gaps of $S$ is equidistributed modulo $m$ if and only if 
\[
    \Ap(S;a)\setminus\{0\}\equiv\{1,2,\dots,a-1\}\pmod*{m},
\]
where $\Ap(S;a)$ denotes the Ap\'ery set of $S$ relative to $a$ (Definition \ref{def:apery}). With this, since $\{1,2,\dots,a-1\}$ is an arithmetic sequence, we then investigate the case where $\Ap(S;a)\setminus\{0\}$ is also an arithmetic sequence; i.e., where
\[
    \Ap(S;a)\setminus\{0\}=\{\beta+\delta,\beta+2\delta,\dots,\beta+(a-1)\delta\}
\] 
for some $\beta\in\Z$ and $\delta\in\N$. We obtain Theorem \ref{thm:main-result}, which gives precise conditions to determine whether the set of gaps of such a numerical semigroup $S$ is equidistributed modulo $m$. It comes down to computing $\gcd(a\delta,m)$ and considering the congruence classes of $a$, $\beta$, and $\delta$ modulo $m$.

There is one family of numerical semigroups $S$ for which the nonzero terms of an Ap\'ery set of $S$ form an arithmetic sequence: numerical semigroups of the form 
\[S=\langle a, ha+d, ha+2d, \dots, ha+(a-1)d\rangle,\]
with $h\ge0$, $a,d\ge1$, and $\gcd(a,d)=1$. This family consists of two well-known subfamilies which correspond to the cases where $h=0$ and $h>0$. When $h=0$, we have $S=\langle a,d\rangle$ (more commonly written $S=\langle a,b\rangle$ with $\gcd(a,b)=1$), a numerical semigroup generated by at most two elements. The study of these numerical semigroups dates back to Sylvester \cite{Sylvester1882}. When $h\ge1$, $S$ is generated by a \emph{generalized arithmetic sequence} of length $a$. Numerical semigroups generated by generalized arithmetic sequences (which may contain fewer than $a$ terms) have been investigated more recently, by Lewin \cite{Lewin1975}, Selmer \cite{Selmer1977}, Ritter \cite{Ritter1998}, and Matthews \cite{Matthews2004}.

Finally, we give specific conditions for which the sets of gaps of numerical semigroups in this family are equidistributed modulo $m$. In general (Corollary~\ref{cor:gen-arith-seq-ED}), the set of gaps of $S$ is equidistributed modulo $m$ if and only if $\gcd(ad,m)=1$ and at least one of four congruences holds. We highlight two special cases here. When $h=0$ (Corollary~\ref{cor:emb-dim-2-ED}), we have $S=\langle a,b\rangle$ and the set of gaps of $S$ is equidistributed modulo $m$ precisely when $\gcd(ab,m)=1$ and at least one of $a$ and $b$ is congruent to 1 modulo $m$. When $h=1$ (Corollary~\ref{cor:pure-arith-ED}), then $S$ is an numerical semigroup generated by a (purely) arithmetic sequence of length $a$, and the set of gaps of $S$ is equidistributed modulo $m$ if and only if $\gcd(ad,m)=1$ and at least one of the following holds: $a\equiv1\pmod*{m}$; or $d\equiv-1\pmod*{m}$.

\subsection{Organization}
This paper is organized as follows. In Section~\ref{sec:prelims}, we present preliminary results about equidistributed multisets, multiset congruence modulo $m$, and connections to congruences of associated generating functions modulo $(x^m-1)\Z[x]$. In Section~\ref{sec:numerical-semigroups}, we briefly describe numerical semigroups. Since the set of gaps of a numerical semigroup is finite, we can apply results from Section \ref{sec:prelims} to find a criterion in terms of an Ap\'ery set of $S$ to determine if the set of gaps of $S$ is equidistributed modulo $m$. 

In Section~\ref{sec:nonzero-apery-arithmetic}, we consider the case where $S$ has an Ap\'ery set for which the nonzero terms form an arithmetic sequence. We completely determine when the set of gaps of $S$ is equidistributed modulo $m$ in Theorem~\ref{thm:main-result}. In Section~\ref{sec:special-cases}, we determine the numerical semigroups $S$ which have a nonzero $a\in S$ for which $\Ap(S;a)\setminus\{0\}$ is an arithmetic sequence. We find that $S$ is in the family of numerical semigroups as described above. We conclude the paper by applying Theorem~\ref{thm:main-result} to this family to get explicit results for when the set of gaps of a numerical semigroup in this family is equidistributed modulo $m$.

\section{Preliminaries}\label{sec:prelims}
In this paper, we will often work with intervals of integers. For $a,b\in\Z$ with $a\le b$ we will let 
\[
    \ints{a,b}:=\{n\in\Z : a\le n\le b\}.
\] 
And since we will need to allow for repetition of elements in Section \ref{sec:nonzero-apery-arithmetic}, we work with multisets, denoting them with curly braces. Multisets have the same properties as sets with one exception: repetition is allowed. For example, $\{0,0,1\}\ne\{0,1\}$. As a result, if $A$ and $B$ are finite multisets of cardinalities $\#(A)$ and $\#(B)$, then the cardinality of their union is $\#\left(A\cup B\right)=\#(A)+\#(B)$.

\subsection{Equidistributed multisets}
Throughout this section, let $A$ and $B$ be finite multiset of integers, and let $m\in\mathbb{N}$.

We begin by defining a function which counts the number of elements of a finite multiset in a particular congruence class. 

\begin{defn}\label{def:n_{r,m}}
    For $r\in\Z$,  let $A_{r,m}:=\{a\in A : a\equiv r\pmod*{m}\}$, the submultiset of $A$ of elements which are congruent to $r$ modulo $m$; and let 
    $
        \nn_{r,m}(A) := \#(A_{r,m}),
    $ 
    the number of elements of $A$ that are congruent to $r$ modulo $m$.
\end{defn}

Since there are $m$ disjoint congruence classes modulo $m$, we have that the multisets $A_{1,m}, A_{2,m}, \dots, A_{m,m}$ form a partition of $A$. Additionally, if $A$ and $B$ are finite multisets, then $\left(A\cup B\right)_{r,m}=\left(A_{r,m}\right)\cup\left(B_{r,m}\right)$, and thus $\nn_{r,m}\left(A\cup B\right)=\nn_{r,m}\left(A\right)+\nn_{r,m}\left(B\right)$. (Put another way, $\nn_{r,m}$ is a monoid homomorphism from the monoid of finite multisets under multiset union to the monoid of nonnegative integers under addition.)

Next, we define multiset congruence modulo $m$.

\begin{defn}\label{def:set-equiv}
    Given finite multisets $A, B\subset\Z$, we say \emph{$A$ and $B$ are congruent multisets modulo $m$}, denoted 
    \[
        A\equiv B\pmod*{m},
    \] 
    if $\nn_{r,m}(A)=\nn_{r,m}(B)$ for all $r$ in an interval of $m$ consecutive integers.
\end{defn}

We will typically take the interval of $m$ consecutive integers to be $\ints{1,m}$ or $\ints{0,m-1}$.

We will now some useful properties of multiset congruence. We'll first show that two congruent multisets must have the same cardinality.

\begin{prop}\label{prop:multiset-congruence-cardinalities-equal}
    If $A\equiv B\pmod*{m}$, then $\#(A)=\#(B)$.
\end{prop}
\begin{proof}
    We assume $A\equiv B\pmod*{m}$, which means $\nn_{r,m}(A)=\nn_{r,m}(B)$ for $r=1,2,\dots,m$. Since $A$ is partitioned by $A_{1,m}, A_{2,m}, \dots, A_{m,m}$, and since $B$ is partitioned by $B_{1,m}, B_{2,m}, \dots, B_{m,m}$, we have 
    \[
    \#(A)
    =\sum\limits_{r=1}^m \nn_{r,m}(A)
    =\sum\limits_{r=1}^m \nn_{r,m}(B)
    =\#(B).
    \]
\end{proof}

Just as with the usual (integer) congruence modulo $m$, multiset congruence modulo $m$ is an equivalence relation. The proof of the following is straightforward.

\begin{prop}\label{prop:equiv reln}
    For finite multisets $A,B,C\subset\mathbb{Z}$ and any $m\in\N$, 
    \begin{itemize}
        \item $A\equiv A\pmod*{m}$;
        \item if $A\equiv B\pmod*{m}$, then $B\equiv A\pmod*{m}$; and
        \item if $A\equiv B\pmod*{m}$ and $B\equiv C\pmod*{m}$, then $A\equiv C\pmod*{m}$.
    \end{itemize}
\end{prop}

Multiset congruence behaves well with respect to multiset unions.

\begin{prop}\label{prop:multiset-union-congruence}
    Let $A,B,C,D$ be finite multisets and suppose that $C\equiv D\pmod*{m}$. Then $A\equiv B\pmod*{m}$ if and only if $A\cup C\equiv B\cup D\pmod*{m}$.
\end{prop}
\begin{proof}
    Observe that $\nn_{r,m}\left(A\cup C\right)-\nn_{r,m}\left(A\right)=\nn_{r,m}\left(C\right)$ and $\nn_{r,m}\left(B\cup D\right)-\nn_{r,m}\left(B\right)=\nn_{r,m}\left(D\right)$. Since $C\equiv D\pmod*{m}$, we have $\nn_{r,m}\left(C\right)=\nn_{r,m}\left(D\right)$ for all $r\in\ints{1,m}$. Hence 
    \[
    \nn_{r,m}\left(A\cup C\right)-\nn_{r,m}\left(A\right) 
    = \nn_{r,m}\left(B\cup D\right)-\nn_{r,m}\left(B\right)
    \]
    for all $r\in\ints{1,m}$.
    Rearranging, 
    \[
    \nn_{r,m}\left(A\cup C\right)-\nn_{r,m}\left(B\cup D\right)
    = \nn_{r,m}\left(A\right) -\nn_{r,m}\left(B\right)
    \]
    for all $r\in\ints{1,m}$.
    This implies that $\nn_{r,m}(A)=\nn_{r,m}(B)$ for all $r\in\ints{1,m}$ if and only if $\nn_{r,m}\left(A\cup C\right)=\nn_{r,m}\left(B\cup D\right)$ for all $r\in\ints{1,m}$. The result follows.
\end{proof}

In particular, if we have two elements $\alpha$ and $\beta$ which are congruent to each other modulo $m$, then we can add them to or remove them from a pair of congruent multisets to obtain a new pair of congruent multisets. More specifically, applying Proposition~\ref{prop:multiset-union-congruence} with $C=\{\alpha\}$ and $D=\{\beta\}$, we obtain the following.

\begin{cor}\label{cor:add-elt-keep-congruence}
    Let $\alpha,\beta\in\Z$ with $\alpha\equiv\beta\pmod*{m}$. For finite multisets $A,B$, we have $A\equiv B\pmod*{m}$ if and only if $A\cup\{\alpha\}\equiv B\cup\{\beta\}\pmod*{m}$.
\end{cor}

We now define an equidistributed multiset modulo $m$.

\begin{defn}\label{defn:equidistributed-modulo-m}
    We say the multiset $A$ is \emph{equidistributed modulo $m$} if $\nn_{r_1,m}(A)=\nn_{r_2,m}(A)$ for all $r_1,r_2\in\ints{1,m}$.
\end{defn}

\begin{lem}\label{lem:m divs card}
    If $A$ is equidistributed modulo $m$, then $m\mid \#(A)$. When this occurs, $\#(A)=m\cdot \nn_{r,m}(A)$ for all $r\in\ints{1,m}$.
\end{lem}
\begin{proof}
    Suppose $A$ is equidistributed modulo $m$. Then $\nn_{r,m}(A)=\nn_{1,m}(A)$ for all $r\in\ints{1,m}$. Since each element of $A$ is in exactly one of the $m$ congruence classes modulo $m$, 
    \[
        \#(A)
        =\nn_{1,m}(A)+\nn_{2,m}(A)+\dots+\nn_{m,m}(A) = m \cdot \nn_{1,m}(A).
    \] 
    Therefore, $m\mid \#(A)$.
\end{proof}

\begin{remark}
As in Definition~\ref{def:n_{r,m}}, we can replace the interval $\ints{1,m}$ in Definition~\ref{defn:equidistributed-modulo-m} and in Lemma~\ref{lem:m divs card} with any interval of $m$ consecutive integers.
\end{remark}

Next, we see that if $A$ is equidistributed modulo $m$, then $A$ is equidistributed modulo all divisors of $m$.

\begin{prop}\label{prop:even-dist-mod-divisor}
    For $d,m\in\N$, if $A$ is equidistributed modulo $m$ and $d\mid m$, then $A$ is equidistributed modulo $d$.
\end{prop}
\begin{proof}
    Suppose $m=dk$ for $k\in\N$. For $n\in\Z$, we have $n\equiv r\pmod*{d}$ if and only if $n\equiv r+id\pmod*{m}$ for some $i\in\ints{1,k}$. In other words, $n\in A_{r,d}$ if and only if $n\in A_{r+id,m}$ for some $i\in\ints{1,k}$. Since $A_{r_1,m}\cap A_{r_2,m}=\emptyset$ for $r_1,r_2\in\ints{1,m}$ with $r_1\ne r_2$,
    \[
        \nn_{r,d}(A)=\sum\limits_{i=1}^k\nn_{r+id,m}(A).
    \]
    Now, if $A$ is equidistributed modulo $m$, then $\nn_{r,m}(A)=\#(A)/m$ for all $r\in\ints{1,m}$. Thus 
    \[
        \nn_{r,d}(A)
        =k\cdot \nn_{r,m}(A)
        =k\cdot\#(A)/m
        =\#(A)/d
    \]
    for all $r\in\ints{1,d}$. We therefore conclude that $A$ is equidistributed modulo $d$.
\end{proof}

\begin{example}\label{ex:ED-example}
    Consider the set $A=\{1,2,3,4,6,8,9,11,13,16,18,23\}$, which has $\#(A)=12$. (This is the set of positive integers which cannot be written as a non-negative linear combination of 5 and 7.) By Lemma~\ref{lem:m divs card}, if $A$ is equidistributed modulo $m$, then $m$ is a divisor of 12. We immediately see that $A$ is not equidistributed modulo 12 because $A$ contains no element which is 0 modulo 12. However, $A$ is equidistributed modulo 6 and modulo 4. By Proposition~\ref{prop:even-dist-mod-divisor}, $A$ is also equidistributed modulo any divisor of 6 or 4. Thus, $A$ is equidistributed modulo $m$ for $m\in\{1,2,3,4,6\}$.
    
    In Example~\ref{ex:ED-example-revisited}, we will revisit this set.
\end{example}

\subsection{Generating functions}
We now consider generating functions arising from multisets. For a finite multiset $A\subset\N_0$, consider the (polynomial) generating function 
\[
    \PP_A(x):=\sum\limits_{a\in A}x^a\in\Z[x].
\] 
For example, if $A=\{3,3,4\}$, then $\PP_A(x)=2x^3+x^4$.

We have thus far worked with multiset congruences in $\Z/m\Z$. The analogue for generating functions is to work with polynomial congruences. For any $h(x)\in\Z[x]$, let $\langle h(x)\rangle:=h(x)\Z[x]$, the principal ideal generated by $h(x)$. (In the introduction, we used angle brackets to denote the set of all finite non-negative linear combinations of a set of integers. For the remainder of this paper, the context will be clear.) As usual, for $f(x),g(x),h(x)\in\Z[x]$, we write 
\[f(x)\equiv g(x)\pmod*{h(x)}\]
to mean $f(x)-g(x)\in\langle h(x)\rangle$.

Now, we fix $m\in\N$. In the polynomial ring $\Z[x]$, consider the ideal $\langle x^m-1\rangle$. We have a ring homomorphism 
\begin{align*} 
\phi : &\ \Z[x]\to\Z[x]/\langle x^m-1\rangle \\
&\ f(x) \mapsto f(x)+\langle x^m-1\rangle.
\end{align*}
Since $(x^m-1)$ is a monic polynomial and $\mathbb{Z}$ is an integral domain, we can apply the polynomial division algorithm to any $f(x)\in\Z[x]$ and $(x^m-1)$ to find a unique pair of polynomials $q_f(x),r_f(x)\in\Z[x]$ for which $f(x)=q_f(x)\cdot(x^m-1)+r_f(x)$ where $\deg(r_f(x))<m$ or $r_f(x)$ is the zero polynomial. In other words, $f(x)+\langle x^m-1\rangle = r_f(x)+\langle x^m-1\rangle$. The uniqueness of the pair means we can represent elements of $\Z[x]/\langle x^m-1\rangle$ with polynomials $r(x)$ of degree less than $m$ along with the zero polynomial:
\[\Z[x]/\langle x^m-1\rangle = \{r(x) + \langle x^m-1\rangle : r(x)\in\Z[x],\, \deg(r(x))<m\}\cup\{0+\langle x^m-1\rangle\}.\]
In particular, since $x^m=1(x^m-1) + 1$, we have $x^m\equiv1\pmod*{(x^m-1)}$.

For this particular ideal, it is straightforward to take a polynomial $f(x)$ of any degree and compute $r_f(x)$. Suppose \[f(x)=\sum\limits_{k=0}^\infty c_k x^k\] for integers $c_0,c_1,\dots$, where all but finitely many are 0. We obtain the image of $f(x)$ in $\Z[x]/\langle x^m-1\rangle$ by repeatedly replacing any instance of $x^m$ with 1. (Put another way, we reduce the exponents modulo $m$.) The result is a polynomial of degree less than $m$. From the previous paragraph, we know there is only one such polynomial, so this reduction must be $r_f(x)$. The coefficient of $x^k$ in $r_f(x)$ is the sum of the coefficients of $x^k, x^{k+m}, x^{k+2m}, \dots$ of $f(x)$. This is the sum of coefficients for which the index is congruent to $k$ modulo $m$.  We have 
\[r_f(x) = \sum\limits_{k=0}^{m-1}\left(\sum\limits_{i=0}^\infty c_{k+im}\right) x^k.\] (All but finitely many of the $c_{k+im}$ terms are zero.) 

Now, suppose $A$ is a finite multiset of nonnegative integers. Then, if we write $\PP_A(x)=\sum\limits_{k=0}^\infty c_k x^k$, we see that for any $k\in\ints{0,m-1}$,
\[\sum\limits_{i=0}^\infty c_{k+im}=c_k+c_{k+m}+c_{k+2m}+\dots=\#(A_{k,m})=\nn_{k,m}(A).\] We immediately obtain the following lemma.

\begin{lem}\label{lem:P_A(x)-as-summation}
If $A$ is a finite multiset of nonnegative integers, then \[\PP_A(x)=\sum\limits_{a\in A}x^a \equiv \sum\limits_{r=0}^{m-1} \nn_{r,m}(A)x^r\pmod*{(x^m-1)}.\]
\end{lem}

We now describe a connection between congruence of multisets modulo $m$ and congruence of generating functions modulo $(x^m-1)$.

\begin{prop}\label{prop:set congruence poly congruence}
    For finite multisets $A$ and $B$ of nonnegative integers, $A\equiv B\pmod*{m}$ if and only if $\PP_{A}(x)\equiv \PP_{B}(x)\pmod*{(x^m-1)}$.
\end{prop}
\begin{proof}
    Suppose $A\equiv B\pmod*{m}$. We equivalently have $\nn_{r,m}(A)=\nn_{r,m}(B)$ for all $r$. Since the polynomials 
    \[\sum\limits_{r=0}^{m-1}\nn_{r,m}(A)x^r \text{ and } \sum\limits_{r=0}^{m-1}\nn_{r,m}(B)x^r\] have degree less than $m$, these polynomials are congruent modulo $(x^m-1)$ precisely when their coefficients are equal, which means $\nn_{r,m}(A)=\nn_{r,m}(B)$ for all $r$. Finally, by Lemma \ref{lem:P_A(x)-as-summation}, these polynomials are congruent, respectively, to $\PP_A(x)$ and $\PP_B(x)$ modulo $(x^m-1)$. Thus, $A\equiv B\pmod*{m}$ if and only if $\PP_A(x)\equiv \PP_B(x)\pmod*{(x^m-1)}$.
\end{proof}

By Lemma~\ref{lem:m divs card}, a finite multiset $A$ is equidistributed modulo $m$ if and only if $\nn_{r,m}(A)=\#(A)/m$ for all $r$. Combined with Lemma~\ref{lem:P_A(x)-as-summation}, we obtain the following.

\begin{cor}\label{cor:A ED iff v1}
    $A$ is equidistributed modulo $m$ if and only if
    \[
        \PP_A(x)\equiv\frac{\#(A)}{m}\sum\limits_{r=0}^{m-1}x^r\pmod*{(x^m-1)}.
    \]
\end{cor}

The summation in the above corollary will come up in a few contexts moving forward, so we give it a name. For any $n\in\N$, let 
\[\CC_n(x):=\dfrac{x^n-1}{x-1}=\sum\limits_{r=0}^{n-1}x^r\in\Z[x].\] 

\begin{prop}\label{prop:A ED mod m poly congruence}
    $A$ is equidistributed modulo $m$ if and only if 
    \[
        (x-1)\PP_A(x)\equiv0\pmod*{(x^m-1)}.
    \]
\end{prop}
\begin{proof}
    For the forward direction, suppose $A$ is equidistributed modulo $m$. Since
    \[
        (x-1)\sum\limits_{r=0}^{m-1}x^r=x^m-1\equiv0\pmod*{(x^m-1)},
    \]
    we use Corollary \ref{cor:A ED iff v1} to conclude that
    \[
        (x-1)\PP_A(x)\equiv(x-1)\dfrac{\#(A)}{m}\sum\limits_{r=0}^{m-1}x^r\equiv 0\pmod*{(x^m-1)},
    \]
    as desired.
    
    For the reverse direction, suppose 
    \[(x-1)\PP_A(x)\equiv0\pmod*{(x^m-1)}.\]
    Thus, $(x-1)\PP_A(x)=(x^m-1)h(x)$ for some $h(x)\in\Z[x]$.
    Dividing through by $(x-1)$, we find $\PP_A(x)=\CC_m(x)h(x)$. We then apply the division algorithm to $h(x)$ and $(x-1)$ (which is valid in $\Z[x]$ because $(x-1)$ is monic) to find 
    \[
        h(x)=(x-1)q(x)+h(1)
    \] 
    for some $q(x)\in\Z[x]$. Thus 
    \[
        \PP_A(x)=\CC_m(x)h(1)+(x^m-1)q(x),
    \] 
    so $\PP_A(x)\equiv \CC_m(x) h(1)\pmod*{(x^m-1)}$. We therefore have $\nn_{r,m}(A)=h(1)$ for all $r\in\ints{1,m}$. Thus, $A$ is equidistributed modulo $m$.
\end{proof}

\section{Numerical semigroups and sets of gaps}\label{sec:numerical-semigroups}
Numerical semigroups were briefly described in the introduction. We provide more details here. For a thorough treatment, see \cite{RosalesGarciaSanchez09}.

A \emph{numerical semigroup} is a submonoid of $\N_0$ under addition with finite complement. In other words, a numerical semigroup is a set $S\subseteq\N_0$ that is closed under addition, contains 0, and has finite complement in $\N_0$. Let $\HH(S)=\mathbb{N}_0\setminus S$, the set of \emph{gaps} of $S$. The \emph{genus} of $S$ is $\genus(S)=\#(\HH(S))$, the number of gaps of $S$. Since $\HH(S)$ is a finite set, 
\[\PP_{\HH(S)}(x)=\sum\limits_{n\in\HH(S)}x^n\in\Z[x].\] 
We call $\PP_{\HH(S)}(x)$ the \emph{gap polynomial of $S$}.

For any set $G\subset\N_0$, let $\langle G\rangle$ denote the set of all finite $\N_0$-linear combinations of elements of $G$. If $G=\{g_1,\dots,g_k\}$, then we let $\langle g_1,\dots,g_k\rangle=\langle G\rangle$. It is a standard result that $\langle G\rangle$ is a numerical semigroup if and only if $\gcd(G)=1$. For any numerical semigroup $S$, there is a finite set $G\subset\N$ such that $S=\langle G\rangle$. Further, there is a unique minimal (relative to set inclusion) generating set $G$. The \emph{embedding dimension} of $S$, denoted $\embdim(S)$, is the cardinality of the minimal generating set $G$ (which is necessarily finite). The least nonzero element of $S$, which is also the least element of its minimal generating set $G$, is the \emph{multiplicity} of $S$, denoted $\mult(S)$. 

An \emph{Ap\'ery set} of a numerical semigroup is an incredibly useful object. We have the following definition.
\begin{defn}[\cite{Apery1946}]\label{def:apery}
    For $S$ a numerical semigroup and any nonzero $a\in S$, the \emph{Ap\'ery set of $S$ relative to $a$} is \[\Ap(S;a)=\left\{s\in S : s-a\not\in S\right\}.\] 
\end{defn}

Put another way, $\Ap(S;a)$ consists of the least element of $S$ in each congruence class modulo $a$. Hence, $\#(\Ap(S;a))=a$. 

\subsection{Detecting an equidistributed set of gaps via an Ap\'ery set}
In the introduction of this paper, we mentioned the work of Wang \& Wang \cite{WangWang2008} from which, for a numerical semigroup of the form $S=\langle a,b\rangle$ with $\gcd(a,b)=1$, one can determine whether or not $\HH(S)$ is equidistributed modulo 2. Their work was based on the following result of Tuenter.
\begin{thm}[{\cite[Theorem 2.1]{Tuenter06}}]\label{thm:tuenter}
    For $S=\langle a,b\rangle$ and any function $f$ defined on $\N_0$, 
    \[
        \sum\limits_{n\in\HH(S)}\left[f(n+a)-f(n)\right]=\sum\limits_{n=1}^{a-1}\left[f(nb)-f(n)\right].
    \]
\end{thm}
Still with $S=\langle a,b\rangle$, Wang \& Wang used Theorem~\ref{thm:tuenter} and the function $f(n)=(-1)^n$ to derive an expression (\cite[Theorem 4.1]{WangWang2008}) for the \emph{alternate Sylvester sum} \[T_m(S)=\sum\limits_{n\in\HH(S)}(-1)^n n^m\] for $m\in\N_0$. In particular, when $m=0$, they found
\[
    \sum\limits_{n\in\HH(S)}(-1)^n = 
    \begin{dcases*}
        0 & for $a,b$ odd,\\
        -\dfrac{b-1}{2} & $a$ even, $b$ odd.
    \end{dcases*}
\]
(Note that $a$ and $b$ are interchangeable. Since they cannot both be even, we may assume that $b$ is odd.) Thus, the numerical semigroup $S=\langle a,b\rangle$ has as many even gaps as odd gaps precisely when $a$ and $b$ are odd. Otherwise, there are more odd gaps than even gaps.

Our goal here is to determine when the set of gaps of a numerical semigroup is equidistributed modulo $m$. We take a similar approach, beginning with a generalization of Tuenter's result. 

\begin{thm}[{\cite[Theorem 2.3]{GassertShor17}}]\label{thm:tuenter gen}
    For $S$ a numerical semigroup, any nonzero $a\in S$, and any function $f$ defined on $\mathbb{N}_0$,
    \[
        \sum\limits_{n\in \HH(S)}[f(n+a)-f(n)] = \sum\limits_{n\in\Ap(S;a)}f(n)-\sum\limits_{n=0}^{a-1}f(n).
    \]
\end{thm}

Since $-1$ is a 2nd root of unity, the function $f(n)=(-1)^n$ is useful for understanding the congruence classes of the gaps modulo 2. In order to understand the congruence classes of the gaps of a numerical semigroup modulo $m$, we should therefore consider $m$th roots of unity. For $\zeta_m$ a primitive $m$th root of unity, we can use $f(n)=\zeta_m^n$. However, since we do not need any particular properties of the cyclotomic ring $\Z[\zeta_m]$, we will stick with polynomials, using $f(n)=x^n$ and working modulo $(x^m-1)$.

With $f(n)=x^n$, Theorem~\ref{thm:tuenter gen} gives the following corollary.
\begin{cor}\label{cor:(x^a-1)P(x)}
    For $S$ a numerical semigroup and any nonzero $a\in S$, 
    \[(x^a-1)\PP_{\HH(S)}(x) = \sum\limits_{n\in\Ap(S;a)}x^n - \sum\limits_{n=0}^{a-1}x^n.\]
\end{cor}

In Proposition~\ref{prop:A ED mod m poly congruence}, we found that a finite multiset $A$ is equidistributed modulo $m$ if and only if \[
    (x-1)\PP_A(x)\equiv0\pmod*{(x^m-1)}.
\]
Since $(x-1)\mid(x^a-1)$ in $\Z[x]$, we can combine that result with Corollary~\ref{cor:(x^a-1)P(x)} to say something about when the set of gaps of $S$ is equidistributed modulo $m$.

\begin{cor}\label{cor:ed-implies-congruence}
    For $S$ a numerical semigroup with some nonzero $a\in S$, if $\HH(S)$ is equidistributed modulo $m$, then 
    \[
        \sum\limits_{n\in\Ap(S;a)}x^n \equiv \sum\limits_{n=0}^{a-1}x^n\pmod*{(x^m-1)},
    \]
    which we may equivalently write as
    $\Ap(S;a)\equiv\ints{0,a-1}\pmod*{m}$.
\end{cor}
\begin{proof}
    Suppose $\HH(S)$ is equidistributed modulo $m$. By Proposition \ref{prop:A ED mod m poly congruence}, 
    \[(x-1)\PP_{\HH(S)}(x)\equiv0\pmod*{(x^m-1)}.\]
    Thus, 
    \[(x^a-1)\PP_{\HH(S)}(x)=\CC_a(x)(x-1)\PP_{\HH(S)}(x)\equiv0\pmod*{(x^m-1)}.\]
    By Corollary \ref{cor:(x^a-1)P(x)}, we conclude that \[\sum\limits_{n\in\Ap(S;a)}x^n\equiv\sum\limits_{n=0}^{a-1}x^n\pmod*{(x^m-1)},\]
    as desired.
    
    By Proposition~\ref{prop:set congruence poly congruence}, this conclusion is equivalent to the statement that $\Ap(S;a)\equiv\ints{0,a-1}\pmod*{m}$.
    \end{proof}

We illustrate this result with an example.

\begin{example}
    Let $S=\langle 3,5\rangle=\{0,3,5,6,8,9,10,\dots\}$, a numerical semigroup with $\HH(S)=\{1,2,4,7\}$. Observe that $\nn_{r,4}(\HH(S))=1$ for each $r\in\ints{0,3}$, which means $\HH(S)$ is equidistributed modulo 4. By Corollary~\ref{cor:ed-implies-congruence}, $\Ap(S;a)\equiv\ints{0,a-1}\pmod*{4}$ for all $a\in S$.
    
    We verify this for a few elements of $S$: 
    \begin{itemize}
        \item $\Ap(S;3)=\{0,5,10\}$;
        \item $\Ap(S,5)=\{0,3,6,9,12\}$; and 
        \item $\Ap(S;14)=\{0,3,5,6,8,9,10,11,12,13,15, 16, 18, 21\}$.
    \end{itemize}
    If we look at each of these multisets modulo 4, we find
    \begin{itemize}
        \item $\Ap(S;3)\equiv\{0,1,2\}\equiv\ints{0,2}\pmod*{4}$,
        \item $\Ap(S;5)\equiv\{0,3,2,1,0\}\equiv\ints{0,4}\pmod*{4}$, and 
        \item $\Ap(S;14)\equiv\{0,3,1,2,0,1,2,3,0,1,3,0,2,1\}\equiv[0,13]\pmod*{4}$.
    \end{itemize}
\end{example}

Next, we have a non-example.

\begin{example}\label{ex:non-example}
    Let $S$ be a numerical semigroup of genus $\genus(S)=g>0$ with $g\in S$. Then $\Ap(S;g)$ is a set of $g$ elements, each in its own congruence class modulo $g$. Hence, $\Ap(S;g)\equiv\ints{0,g-1}\pmod*{g}$. By Lemma~\ref{lem:m divs card}, if $\HH(S)$ is equidistributed modulo $g$, then 
    \[
        \nn_{0,g}(\HH(S))=\dfrac{\#(\HH(S))}{g}=\dfrac{g}{g}=1.
    \] 
    However, since $0\cdot g, 1\cdot g, 2\cdot g, \dots$ are all in $S$, $\nn_{0,g}(\HH(S))=0$. Thus, 
    $\HH(S)$ is not equidistributed modulo $g$.
    
    For a concrete non-example, let $S=\langle 4,5,11\rangle$. The set of gaps of $S$ is $\HH(S)=\{1,2,3,6,7\}$ and $\genus(S)=5\in S$. Observe that $\Ap(S;5)=\{0,4,8,11,12\}$, which has $\Ap(S;5)\equiv\ints{0,4}\pmod*{5}$, yet $\HH(S)$ is not equidistributed modulo 5.
\end{example}

Thus, the converse of Corollary~\ref{cor:ed-implies-congruence} is false. However, if we further include the assumption that $\gcd(a,m)=1$, then the converse holds. We will prove it after stating the following lemma, which is well-known.

\begin{lem}\label{lem:poly-gcd}
    Let $a,b\in\N$. Then $\gcd(x^a-1,x^b-1)=x^{\gcd(a,b)}-1$.
\end{lem}

\begin{prop}\label{prop:ED iff congruence}
    For $S$ a numerical semigroup, any nonzero $a\in S$, and any $m\in\mathbb{N}$ with $\gcd(a,m)=1$, the following statements are equivalent.
    \begin{enumerate} 
        \item $\HH(S)$ is equidistributed modulo $m$.
        \item $\PP_{\HH(S)}(x)\in\langle C_m(x)\rangle$.
        \item $\displaystyle
        \sum\limits_{n\in\Ap(S;a)}x^n \equiv \sum\limits_{n=0}^{a-1}x^n\pmod*{(x^m-1)}$.
        \item $\displaystyle
        \sum\limits_{n\in\Ap(S;a)\setminus\{0\}}x^n \equiv \sum\limits_{n=1}^{a-1}x^n\pmod*{(x^m-1)}$.
        \item $\displaystyle\Ap(S;a)\equiv\ints{0,a-1}\pmod*{m}$.
        \item $\displaystyle\Ap(S;a)\setminus\{0\}\equiv\ints{1,a-1}\pmod*{m}$.
    \end{enumerate}
\end{prop}
\begin{proof}
    By Proposition~\ref{prop:A ED mod m poly congruence}, (1)$\iff$(2). Since $0\in\Ap(S;a)$, (3)$\iff$(4) and (5)$\iff$(6). By Proposition~\ref{prop:set congruence poly congruence}, (3)$\iff$(5). By Corollary~\ref{cor:ed-implies-congruence}, (1)$\implies$(3). It is therefore enough to show (3)$\implies$(1).
    
    Since $\gcd(a,m)=1$, by Lemma~\ref{lem:poly-gcd},  $\gcd(x^a-1,x^m-1)=x^{1}-1$. Therefore, $\CC_a(x)$ is a unit modulo $(x^m-1)$.
    
    Now, suppose
    \[
        \sum\limits_{n\in\Ap(S;a)}x^n\equiv\sum\limits_{n=0}^{a-1}x^n\pmod*{(x^m-1)}.
    \] 
    By Corollary \ref{cor:(x^a-1)P(x)}, $(x^a-1)\PP_{\HH(S)}(x)\equiv0\pmod*{(x^m-1)}$. Therefore 
    \[
        \CC_a(x)(x-1)\PP_{\HH(S)}(x)\equiv0\pmod*{(x^m-1)}.
    \]
    Since $\CC_a(x)$ is a unit modulo $(x^m-1)$, we conclude that $(x-1)\PP_{\HH(S)}(x)\equiv0\pmod*{(x^m-1)}$. By Proposition~\ref{prop:A ED mod m poly congruence}, we conclude that $\HH(S)$ is equidistributed modulo $m$.
\end{proof}

This leads to the following interesting result.

\begin{cor}
    For any nonzero $a,b\in S$ with $\gcd(ab,m)=1$, $\Ap(S;a)\equiv\ints{0,a-1}\pmod*{m}$ if and only if $\Ap(S;b)\equiv\ints{0,b-1}\pmod*{m}$.
\end{cor}

Before moving on, we make one more remark regarding $m$th roots of unity. By Proposition~\ref{prop:A ED mod m poly congruence}, $\HH(S)$ is equidistributed modulo $m$ if and only if $(x-1)\PP_{\HH(S)}(x)\equiv0\pmod*{(x^m-1)}$, which is to say $\CC_m(x)$ divides the gap polynomial $\PP_{\HH(S)}(x)$. For $\zeta_m$ a primitive $m$th root of unity, 
\[
    \CC_m(x)=\prod\limits_{i=1}^{m-1}\left(x-\zeta_m^i\right).
\]
Thus, $\HH(S)$ is equidistributed modulo $m$ if and only if $\PP_{\HH(S)}\left(\zeta_m^i\right)=0$ for all $i=1,2,\dots,m-1$.

\subsection{Numerical semigroups with maximal embedding dimension}

We turn our attention to numerical semigroups of maximal embedding dimension.

Let $S$ be a numerical semigroup. Then $S$ is generated by a minimal set $A$ of cardinality $\embdim(S)$. By \cite[Proposition 2.10]{RosalesGarciaSanchez09},  $\embdim(S)\le\mult(S)$. If $\embdim(S)=\mult(S)$, then $S$ has \emph{maximal embedding dimension}. Let $b\in A$ and $s\in S\setminus\{0,b\}$. By the minimality of $A$, $b-s\not\in A$ and therefore $b\in\Ap(S;s)$. In particular, if $a\in A$,  then $(A\setminus\{a\})\subset\Ap(S;a)$. For $S$ with maximal embedding dimension, we can describe $\Ap(S;a)$ in the case where $a=\mult(S)$.

\begin{prop}[{\cite[Proposition 3.1]{RosalesGarciaSanchez09}}]\label{prop:max-emb-dim-apery-set}
    For $S$ generated by the minimal set $A$, let $a=\mult(S)$. Then $S$ has maximal embedding dimension if and only if 
    \[\Ap(S;a)\setminus\{0\}=A\setminus\{a\}.\]
\end{prop}

We combine Proposition \ref{prop:max-emb-dim-apery-set} with Proposition \ref{prop:ED iff congruence} for a simple criterion to determine when the gaps of a numerical semigroup of maximal embedding dimension are equidistributed modulo $m$.

\begin{cor}\label{cor:max-emb-dim-gaps-ed-iff}
    Suppose $S$ is a numerical semigroup of maximal embedding dimension with minimal generating set $A$. For $a=\mult(S)=\embdim(S)$ and $m\in\N$ with $\gcd(a,m)=1$, $\HH(S)$ is equidistributed modulo $m$ if and only if 
    \[A\setminus\{a\}\equiv\ints{1,a-1}\pmod*{m}.\]
\end{cor}

\subsection{Numerical semigroups with small multiplicity}
We conclude this section by completely determining when the set of gaps of a numerical semigroup of multiplicity 2 or 3 is equidistributed modulo $m$.

To begin, suppose $\mult(S)=2$. Then $\genus(S)>0$, $1\not\in S$, and $2\in S$. Therefore, $\embdim(S)=2$, which means $S$ is a semigroup of maximal embedding dimension. In particular, $S=\langle 2,b\rangle$ for odd $b\ge3$. We can then determine when $\HH(S)$ is equidistributed modulo $m$.

\begin{prop}\label{prop:S-contains-2}
    Suppose $\mult(S)=2$. Then $S=\langle 2,b\rangle$ for some odd integer $b\ge3$, and $\HH(S)$ is equidistributed modulo $m$ if and only if $m$ is an odd divisor of $\genus(S)=(b-1)/2$.
\end{prop}
\begin{proof}
    Suppose $\mult(S)=2$. Then $\embdim(S)=2$. We have $S=\langle 2,b\rangle$ for some $b>1$ with $\gcd(b,2)=1$. In other words, $b$ is an odd integer with $b\ge3$. The set of gaps of $S$ is $\HH(S)=\{1,3,5,\dots,b-2\}$, a set of cardinality $(b-1)/2$. (Thus, $\genus(S)=(b-1)/2$.)
    
    Now, we wish to determine when $\HH(S)$ is equidistributed modulo $m$. Since $\#(\HH(S))=(b-1)/2$, it is a necessary condition (by Lemma \ref{lem:m divs card}) that $m\mid(b-1)/2$. We consider the cases of $m$ even and $m$ odd separately.
    
    First, consider $\HH(S)$ modulo 2. Since $\HH(S)$ contains zero even elements, \[\nn_{0,2}(\HH(S))=0\ne(b-1)/2=\nn_{1,2}(\HH(S)).\]
    Therefore, $\HH(S)$ is not equidistributed modulo 2.
    By the contrapositive of  Proposition~\ref{prop:even-dist-mod-divisor}, $\HH(S)$ is not equidistributed modulo $m$ for $m$ even.
    
    Next, factor out all powers of 2 from $(b-1)/2$ to write
    $(b-1)/2=2^k l$
    for $k\in\N_0$ and $l$ odd. Then $b=1+2^{k+1}l$ and $\gcd(2,l)=1$. By Corollary~\ref{cor:max-emb-dim-gaps-ed-iff}, since $A\setminus\{2\}=\{b\}\equiv\{1\}\pmod*{l}$, $\HH(S)$ is equidistributed modulo $l$. If $m$ is any odd divisor of $(b-1)/2$, then $m$ is a divisor of $l$. Therefore, by Proposition~\ref{prop:even-dist-mod-divisor}, $\HH(S)$ is equidistributed modulo $m$.
    
    We conclude that $\HH(S)$ is equidistributed modulo $m$ if and only if $m$ is an odd divisor of $\#(\HH(S))=(b-1)/2$, as desired.
\end{proof}

\begin{remark}\label{rmk:mult-emb-dim-2}
We note that the previous proposition is a special case of a more general result concerning the equidistribution of sets of gaps of numerical semigroups of embedding dimension 2. (See Corollary \ref{cor:emb-dim-2-ED}.)
\end{remark} 

Next, we consider the case where $\mult(S)=3$.

\begin{prop}\label{prop:S-contains-3}
    Suppose $\mult(S)=3$. Then $\embdim(S)=2$ or 3. 
    
    If $\embdim(S)=2$, then $S=\langle 3,b\rangle$ for some $b\ge3$ with $\gcd(3,b)=1$, and $\HH(S)$ is equidistributed modulo $m$ if and only if $\gcd(3,m)=1$ and $m$ is a divisor of $b-1$.
    
    If $\embdim(S)=3$, then $S=\langle 3,b,c\rangle$ for some $b$, $c$ with $\gcd(3,bc)=1$ and $3<b<c<2b$, and $\HH(S)$ is equidistributed modulo $m$ if and only if $\gcd(3,m)=1$ and $m$ is a divisor of either $\gcd(b-1,c-2)$ or $\gcd(b-2,c-1)$.
\end{prop}
\begin{proof}
     Suppose $\HH(S)$ is equidistributed modulo $m$. Since $1\not\in S$, $\embdim(S)>1$, and since $\embdim(S)\le\mult(S)=3$, we conclude that $\embdim(S)=2$ or 3. Since $3\in S$, $\HH(S)$ contains no multiples of 3. Thus $\HH(S)$ is not equidistributed modulo $m$ for $m$ a multiple of 3. We therefore have $\gcd(3,m)=1$.
     
     If $\embdim(S)=2$, then $S=\langle 3,b\rangle$ for $b>3$ and $\gcd(b,3)=1$. Note that $\Ap(S;3)=\{0,b,2b\}$. (This is straightforward to show. Or see Proposition~\ref{prop:emb-dim-2-apery-set}.) By Proposition~\ref{prop:ED iff congruence}, $\HH(S)$ is equidistributed modulo $m$ if and only if $\{b,2b\}\equiv\{1,2\}\pmod*{m}$. We have two possibilities. If $b\equiv1\pmod*{m}$ then $2b\equiv2\pmod*{m}$ and we are done. If $b\equiv2\pmod*{m}$ and $2b\equiv1\pmod*{m}$, then $4\equiv1\pmod*{m}$ and therefore $m\mid 3$. Since $\gcd(m,3)=1$, we conclude $m=1$, which implies $b\equiv1\pmod*{m}$.
     
     If $\embdim(S)=3$, then $S=\langle 3,b,c\rangle$ for $3<b<c$ and $\gcd(3,bc)=1$. We also have $b\not\equiv c\pmod*{3}$. Since there are only two nonzero congruence classes modulo 3, $2b\equiv c\pmod*{3}$, and thus we must have $c<2b$. As this is a maximal embedding numerical semigroup, by Corollary~\ref{cor:max-emb-dim-gaps-ed-iff} $\HH(S)$ is equidistributed modulo $m$ if and only if $\{b,c\}\equiv\{1,2\}\pmod*{m}$. We must have one of the following: $b\equiv1\pmod*{m}$ and $c\equiv2\pmod*{m}$; or $b\equiv2\pmod*{m}$ and $c\equiv1\pmod*{m}$. For the first case, we equivalently have that $m$ divides $b-1$ and $c-2$, which occurs precisely when $m$ divides $\gcd(b-1,c-2)$. The second case gives the other condition.
\end{proof}

\begin{remark}\label{rmk:mult-emb-dim-3}
Note that if $\mult(S)=\embdim(S)=3$, then $S=\langle 3,b,c\rangle$ for $3<b<c<2b$, $\gcd(3,bc)=1$, and $b\not\equiv c\pmod*{m}$. Then $2b-c=3h$ for some $h\in\N$. Let $d=c-b\in\N$. Then $S=\langle 3, 3h+d, 3h+2d\rangle$. In other words, $S$ is a numerical semigroup of maximal embedding dimension generated by a generalized arithmetic sequence with $\mult(S)=3$. In Section~\ref{sec:special-cases}, we will find explicit results for such semigroups of any multiplicity. (See Corollary~\ref{cor:gen-arith-seq-ED}.) 
\end{remark}

\section{Ap\'ery sets where the nonzero terms form an arithmetic sequence}\label{sec:nonzero-apery-arithmetic}
By Proposition~\ref{prop:ED iff congruence}, for $S$ a numerical semigroup with somen nonzero $a\in S$, $\HH(S)$ is equidistributed modulo $m$ if and only if $\Ap(S;a)\setminus\{0\}\equiv\ints{1,a-1}\pmod*{m}$. Since $\ints{1,a-1}$ is an arithmetic sequence, it seems reasonable to assume that we can get explicit results when $\Ap(S;a)\setminus\{0\}$ is an arithmetic sequence. In this section, we will do just that. In the following section, we will describe the numerical semigroups for which this occurs.

Let $S$ be a numerical semigroup, let $a\in S$, and suppose $\Ap(S;a)\setminus\{0\}$ forms an arithmetic sequence with common difference $\delta$. If $\delta<0$, then we may reverse the sequence to obtain a sequence with the same terms and common difference $-\delta>0$. Thus, we may assume the common difference is positive. In other words, we suppose
\begin{equation}\label{eqn:non-zero-apery-notation}
    \Ap(S;a) = \{0\}\cup\{\beta+\delta, \beta+2\delta, \dots, \beta+(a-1)\delta\}
\end{equation}
for some $\beta\in\Z$ and $\delta\in\N$. Our goal is to determine conditions on $a$, $\beta$, $\delta$, and $m$ for which $\HH(S)$ is equidistributed modulo $m$. 

If $a=1$, then $\Ap(S;a)=\{0\}$. We are free to choose $\beta$ and $\delta$ as we wish. In this case, we will let $\beta=0$ and $\delta=1$. If $a=2$, then $\Ap(S;a)=\{0,\beta+\delta\}$. We must have $\beta+\delta$ odd. In this case, we will let $\beta=0$ and $\delta$ be an odd positive integer. For $a\ge3$, the values of $\beta$ and $\delta$ are uniquely determined by the sequence $\left\{\beta+i\delta : i\in\ints{1,a-1}\right\}$, which contains at least two elements.

\begin{lem}\label{lem:arith-props}
    If $\Ap(S;a)=\{0\}\cup\left\{\beta+i\delta: i\in\ints{1,a-1}\right\}$, then $\gcd(a,\delta)=1$, $a\mid \beta$, and $\beta\ge0$.
\end{lem}
\begin{proof}
    The conclusions hold for our choices of $\beta$ and $\delta$ for $a=1$ and $a=2$. For the rest of this proof, we suppose $a\ge3$. Since $\Ap(S;a)\setminus\{0\}\equiv\ints{1,a-1}\pmod*{a}$ and $a\ge3$, there are integers $i,j$ with $1\le i,j\le a-1$ for which $\beta+i\delta\equiv1\pmod*{a}$ and $\beta+j\delta\equiv2\pmod*{a}$. Taking the difference, we have $(j-i)\delta\equiv1\pmod*{a}$, and thus $\gcd(a,\delta)=1$.
    
    Next, since $\gcd(a,\delta)=1$, $\delta$ is a generator of $\mathbb{Z}/a\mathbb{Z}$. In particular, $\{ i\delta : 1\le i\le a-1\}\equiv\ints{1,a-1}\pmod*{a}$. Therefore
    \[
    \ints{\beta+1,\beta+a-1}
    \equiv\{\beta + i\delta : 1\le i\le a-1\}
    \equiv\Ap(S;a)\setminus\{0\}
    \equiv\ints{1,a-1}\pmod*{a},
    \]
    implying $\beta\equiv0\pmod*{a}$ and hence $a\mid\beta$.
    
    Finally, since $\beta+\delta\in S$, $2(\beta+\delta)\in S$ as well. Also, $\beta+2\delta\in\Ap(S;a)\subset S$. We have $a\mid \beta$, and thus $\beta+2\delta\equiv 2(\beta+\delta)\equiv 2\delta\pmod*{a}$. Since $\beta+2\delta\in\Ap(S;a)$, $\beta+2\delta$ is the smallest element of $S$ in its congruence class modulo $a$. Thus, $\beta+2\delta\le 2(\beta+\delta)$, from which we conclude $0\le \beta$.
\end{proof}

Now, if $\Ap(S;a)\setminus\{0\}$ is an arithmetic sequence, then by Corollary \ref{cor:(x^a-1)P(x)}, we have
\begin{equation}\label{eqn:first-version}
    (x^a-1)P_{\HH(S)}(x) 
    = x^\beta\sum\limits_{j=1}^{a-1}x^{j\delta} - \sum\limits_{j=1}^{a-1}x^j.
\end{equation}
Each summation is a finite geometric series, so we may equivalently write 
\begin{equation}
    (x^a-1)P_{\HH(S)}(x) =x^\beta\dfrac{x^{a\delta}-x^\delta}{x^\delta-1}-\dfrac{x^a-x^1}{x-1},
\end{equation}
which is equivalent to the polynomial equation
\begin{equation}\label{eqn:third-version}
    (x-1)(x^a-1)(x^\delta-1)P_{\HH(S)}(x) = x^\beta(x-1)(x^{a\delta}-x^\delta)-(x^a-x)(x^\delta-1).
\end{equation}

We will use these equations to show that if $\HH(S)$ is equidistributed modulo $m$, then $\gcd(\delta,m)=1$ and $\gcd(a,m)=1$.

\begin{prop}\label{prop:gcd(delta,m)=1}
    Suppose $\HH(S)$ is equidistributed modulo $m$. Then $\gcd(\delta,m)=1$.
\end{prop}

\begin{proof}
    If $a=1$, then $\delta=1$ and hence $\gcd(\delta,m)=1$.
    
    For the rest of this proof, we assume $a\ge2$. By Proposition~\ref{prop:A ED mod m poly congruence}, if $\HH(S)$ is equidistributed modulo $m$, then $(x-1)P_{\HH(S)}(x)\equiv0\pmod*{(x^m-1)}$. Thus, by Eq.~\eqref{eqn:first-version}, 
    \[
        x^\beta(x^\delta+x^{2\delta}+\dots+x^{(a-1)\delta}) \equiv x^1+x^2+\dots+x^{a-1}\pmod*{(x^m-1)}.
    \]
    Equivalently, 
    \[\{\beta+\delta,\beta+2\delta,\dots,\beta+(a-1)\delta\}\equiv\{1,2,3\dots,a-1\}\pmod*{m}.\]
    
    If $a=2$, then $\beta+\delta\equiv1\pmod*{m}$. Since $\beta=0$, this implies $\gcd(\delta,m)=1$.
    
    If $a\ge3$, we must have some $i,j\in\ints{1,a-1}$ such that $\beta+i\delta\equiv1\pmod*{m}$ and $\beta+j\delta\equiv2\pmod*{m}$. In particular, $(j-i)\delta\equiv1\pmod*{m}$, from which we conclude that $\gcd(\delta,m)=1$.
\end{proof}

\begin{prop}\label{prop:gcd(a,m)=1}
    Suppose $\HH(S)$ is equidistributed modulo $m$. Then $\gcd(a,m)=1$.    
\end{prop}
\begin{proof}
    Suppose $\HH(S)$ is equidistributed modulo $m$. By Lemma~\ref{lem:arith-props}, $a\mid \beta$. Thus, $\beta=al$ for some $l\in\Z$. By the division algorithm in $\Z$, there exist $q,r\in\Z$ for which $l=qm+r$ with $1<r\le m+1$. (It will be helpful to have $r\ge2$.) It follows that $\beta\equiv ar\pmod*{m}$.
    
    By Proposition \ref{prop:gcd(delta,m)=1}, $\gcd(\delta,m)=1$. Thus, $(x-1)P_{\HH(S)}(x)\equiv0\pmod*{(x^m-1)}$ if and only if 
    \[(x-1)(x^\delta-1)P_{\HH(S)}(x)\equiv0\pmod*{(x^m-1)}.\]
    Since $\HH(S)$ is equidistributed modulo $m$, by Eq.~\eqref{eqn:third-version}, 
    \[\dfrac{x^\beta\left(x-1\right)\left(x^{a\delta}-x^\delta\right)-\left(x^a-x\right)\left(x^\delta-1\right)}{x^a-1}\equiv0\pmod*{(x^m-1)}.\]
    Since $\beta\equiv ar\pmod*{m}$, we expand and regroup to find
    \[
        \dfrac{x\left(x^{(r+\delta)a}-1\right)}{x^a-1}
        +\dfrac{x^{a+\delta}\left(x^{(r-1)a}-1\right)}{x^a-1}
        \equiv
        \dfrac{x^a\left(x^{(r+\delta-1)a}-1\right)}{x^a-1}
        +\dfrac{x^{\delta+1}\left(x^{ra}-1\right)}{x^a-1}
    \pmod*{(x^m-1)}.
    \]
    Each quotient is a finite geometric series. (Note that $r+\delta$, $r-1$, $r+\delta-1$, and $r$ are all positive integers. This is why we wanted $r\ge2$ earlier.) The left side expands out as
    \[
    L(x)=x\left(1+x^a+x^{2a}+\dots+x^{(r+\delta-1)a}\right)
    + x^{a+\delta}\left(1+x^a+x^{2a}+\dots+x^{(r-2)a}\right)
    \]
    and the right side expands out as
    \[
    R(x)=x^a\left(1+x^a+x^{2a}+\dots+x^{(r+\delta-2)a}\right)
    + x^{\delta+1}\left(1+x^a+x^{2a}+\dots+x^{(r-1)a}\right).
    \]
    
    Now, since $L(x)\equiv R(x)\pmod*{(x^m-1)}$, for the corresponding multisets 
    \[
        M_L=\left\{1+ia : 0\le i\le r+\delta-1\right\}\cup\left\{\delta+ja : 1\le j\le r-1\right\}
    \]
    and
    \[
        M_R=\left\{ia : 1\le i\le r+\delta-1\right\}\cup\left\{\delta+1+ja : 0\le j\le r-1\right\},
    \]
    we have $M_L\equiv M_R\pmod*{m}$.
    In particular, since $1\in M_L$, some element of $M_R$ is congruent to $1$ modulo $m$. Thus, we have two cases: either $ia\equiv1\pmod*{m}$ for some $i\in\ints{1,r+\delta-1}$; or $\delta+1+ja\equiv1\pmod*{m}$ for some $j\in\ints{0,r-1}$.
    
    In the first case, we have that $a$ is a unit modulo $m$ and hence $\gcd(a,m)=1$, as desired.
    
    In the second case, we have two subcases. If $j=0$, then $\delta+1\equiv1\pmod*{m}$, which implies $\delta\mid m$. Since $\gcd(\delta,m)=1$, this implies $m=1$, and thus $\gcd(a,m)=1$. If $j>0$, then $\delta+1+ja\equiv1\pmod*{m}$. Since $\gcd(\delta,m)=1$, $\delta$ is a unit modulo $m$ and thus $-\delta^{-1}ja\equiv1\pmod*{m}$. Therefore, $a$ is also a unit modulo $m$, and hence $\gcd(a,m)=1$.
\end{proof}

In the above proof, $M_L$ and $M_R$ may have repeated elements. For this reason, as well as for similar reasons in subsequent results in this section, we have chosen to work with multisets in this paper rather than sets.

\begin{prop}\label{prop:gaps ED iff multiset cong}
    $\HH(S)$ is equidistributed modulo $m$ if and only if $\gcd(a\delta,m)=1$ and we have the multiset congruence 
    \[\{a,a\delta+\beta+1,\beta+\delta,\delta+1\}\equiv\{a+\delta,a\delta+\beta,\beta+\delta+1,1\}\pmod*{m}.\] 
\end{prop}
\begin{proof}
    For the forward direction, suppose $\HH(S)$ is equidistributed modulo $m$. By Proposition \ref{prop:gcd(delta,m)=1} and Proposition \ref{prop:gcd(a,m)=1}, $\gcd(\delta,m)=\gcd(a,m)=1$. By Proposition \ref{prop:A ED mod m poly congruence}, $(x-1)\PP_{\HH(S)}(x)\equiv0\pmod*{(x^m-1)}$, and hence 
    \[
        (x-1)(x^a-1)(x^\delta-1)\PP_{\HH(S)}(x)\equiv0\pmod*{(x^m-1)}.
    \]
    By Eq.\ \eqref{eqn:third-version}, 
    \[
        x^\beta(x-1)(x^{a\delta}-x^\delta)-(x^a-x)(x^\delta-1)\equiv0\pmod*{(x^m-1)}.
    \]
    Expanding out and moving terms so that all leading coefficients are $+1$, we have
    \[
        x^a+x^{a\delta+\beta+1}+x^{\beta+\delta}+x^{\delta+1}\equiv
        x^{a+\delta}+x^{a\delta+\beta}+x^{\beta+\delta+1}+x^1\pmod*{(x^m-1)}.
    \]
    By Proposition~\ref{prop:set congruence poly congruence}, 
    \[
        \{a,a\delta+\beta+1,\beta+\delta,\delta+1\}\equiv
        \{a+\delta,a\delta+\beta,\beta+\delta+1,1\}\pmod*{m},
    \]
    as desired.
    
    Now we prove the converse. Suppose 
    \[
        \{a,a\delta+\beta+1,\beta+\delta,\delta+1\}\equiv
        \{a+\delta,a\delta+\beta,\beta+\delta+1,1\}\pmod*{m}
    \] 
    and $\gcd(a\delta,m)=1$. Reversing our steps from the proof of the forward implication, the multiset congruence implies
    \[
        (x-1)(x^a-1)(x^\delta-1)\PP_{\HH(S)}(x)\equiv0\pmod*{(x^m-1)}.
    \]
    Since $\gcd(a\delta,m)=1$, we use Lemma~\ref{lem:poly-gcd} to get 
    \[
        \gcd((x-1)(x^a-1)(x^\delta-1),x^m-1)=\gcd((x-1)^3,x^m-1)=x-1.
    \]
    We therefore have that 
    \[
        (x-1)\PP_{\HH(S)}(x)\equiv0\pmod*{(x^m-1)},
    \]
    which, with Proposition \ref{prop:A ED mod m poly congruence}, implies that $\HH(S)$ is equidistributed modulo $m$, as desired.
\end{proof}

We are now ready to state the main result. The proof follows from our work above along with tracking through specific cases.

\begin{thm}\label{thm:main-result}
    For $S$ a numerical semigroup, any nonzero $a\in S$, and  
    \[
        \Ap(S;a)=\{0\}\cup\{\beta+i\delta:i\in\ints{1,a-1}\}
    \] 
    for some $\beta\ge0$ and $\delta>0$,
    $\HH(S)$ is equidistributed modulo $m$ if and only if $\gcd(a\delta,m)=1$ and one of the following occurs:
    \begin{enumerate}
        \item $a\equiv1\pmod*{m}$; or
        \item $a\equiv2\pmod*{m}$ and $\beta+\delta\equiv1\pmod*{m}$; or
        \item $\delta\equiv1\pmod*{m}$ and $\beta\equiv0\pmod*{m}$; or
        \item $\delta\equiv-1\pmod*{m}$ and $\beta\equiv a\pmod*{m}$.
    \end{enumerate}
    (Recall that $\beta$ and $\delta$ are uniquely determined when $a\ge3$. For $a=1$, we take $\beta=0$ and $\delta=1$. For $a=2$, we take $\beta=0$ and $\delta$ odd.)
\end{thm}
\begin{proof}
    We begin with the forward implication. If $\HH(S)$ is equidistributed modulo $m$, then by Proposition \ref{prop:gaps ED iff multiset cong}, $\gcd(a\delta,m)=1$ and 
    \[
        \{a,a\delta+\beta+1,\beta+\delta,\delta+1\}\equiv
        \{a+\delta,a\delta+\beta,\beta+\delta+1,1\}\pmod*{m}.
    \]
    In particular, $a$ is congruent modulo $m$ to one of the following elements: $a+\delta$; $a\delta+\beta$; $\beta+\delta+1$; or $1$. We consider these four cases separately.
    
    \begin{enumerate} 
        \item 
        If $a\equiv a+\delta\pmod*{m}$, then $\delta\equiv0\pmod*{m}$. Since $\gcd(\delta,m)=1$, we have $m=1$. All four conclusions therefore hold.
        
        \item 
        If $a\equiv a\delta+\beta\pmod*{m}$, then $\beta\equiv a-a\delta\pmod*{m}$. We use Corollary~\ref{cor:add-elt-keep-congruence}, looking at the remaining elements in the multisets to find 
        $\{a+1,a-a\delta+\delta,\delta+1\}\equiv\{a+\delta,a-a\delta+\delta+1,1\}\pmod*{m}$. Next, $a+1$ is congruent modulo $m$ to one of the following elements: $a-a\delta+\delta+1$; $a+\delta$; or $1$. We consider these cases separately.
        \begin{enumerate}
            \item If $a+1\equiv a+\delta\pmod*{m}$, then $\delta\equiv1\pmod*{m}$, which then implies $\beta\equiv0\pmod*{m}$. This is the third concluding case.
            \item If $a+1\equiv a-a\delta+\delta+1\pmod*{m}$, then $a\delta\equiv \delta\pmod*{m}$. Since $\gcd(\delta,m)=1$, we get $a\equiv1\pmod*{m}$, which is the first concluding case.
            \item If $a+1\equiv1\pmod*{m}$, then $m\mid a$. Since $\gcd(a,m)=1$, we have $m=1$. All four concluding cases occur.
        \end{enumerate}
        
        \item 
        If $a\equiv \beta+\delta+1\pmod*{m}$, we have $\beta\equiv a-\delta-1\pmod*{m}$. We use Corollary ~\ref{cor:add-elt-keep-congruence}, looking at the remaining elements in the multisets to find
        $\{a+a\delta-\delta,a-1,\delta+1\}\equiv\{a+\delta,a+a\delta-\delta-1,1\}\pmod*{m}$. We see that $a-1$ is congruent modulo $m$ to one of $a+\delta, a+a\delta-\delta-1,  1$. We consider these cases separately. 
        \begin{enumerate} 
            \item If $a-1\equiv a+\delta\pmod*{m}$, then $\delta\equiv-1\pmod*{m}$. Since $\beta\equiv a-\delta-1\pmod*{m}$, we have $\beta\equiv a\pmod*{m}$. This is the fourth concluding case.
            \item If $a-1\equiv a+a\delta-\delta-1\pmod*{m}$, then $a\delta\equiv \delta\pmod*{m}$ and hence $a\equiv1\pmod*{m}$. This is the first concluding case.
            \item If $a-1\equiv 1\pmod*{m}$, then $a\equiv2\pmod*{m}$. Since $\beta\equiv a-\delta-1\pmod*{m}$, we have $\beta+\delta\equiv1\pmod*{m}$. This is the second concluding case.
        \end{enumerate}
        \item 
            If $a\equiv1\pmod*{m}$, then we are done. This is the first concluding case.
    \end{enumerate}

    Now we prove the converse. For each of the four concluding cases, we need to show 
    \[
        \{a,a\delta+\beta+1,\beta+\delta,\delta+1\}\equiv
        \{a+\delta,a\delta+\beta,\beta+\delta+1,1\}\pmod*{m}.
    \] 
    Along with the condition that $\gcd(a\delta,m)=1$, we will use Proposition \ref{prop:gaps ED iff multiset cong} to conclude that $\HH(S)$ is equidistributed modulo $m$.
    
    For the following, let $L=\{a,a\delta+\beta+1,\beta+\delta,\delta+1\}$ and $R=\{a+\delta,a\delta+\beta,\beta+\delta+1,1\}$. It is therefore enough to verify that $L\equiv R\pmod*{m}$ in each case.
    
    \begin{enumerate}
    \item If $a\equiv1\pmod*{m}$, then \[
        L\equiv\{1,\delta+\beta+1,\beta+\delta,\delta+1\}\pmod*{m}
    \]
    and 
    \[
        R\equiv\{1+\delta,\delta+\beta,\beta+\delta+1,1\}\pmod*{m}.
    \]
    
    \item If $a\equiv2\pmod*{m}$ and $\beta+\delta\equiv1\pmod*{m}$, then 
    \[
        L\equiv\{2,2\delta+1-\delta+1,1-\delta+\delta,\delta+1\}\equiv \{2,\delta+2,1,\delta+1\}\pmod*{m}
    \] 
    and 
    \[
        R\equiv\{2+\delta,2\delta+1-\delta,1-\delta+\delta+1,1\}\equiv\{2+\delta,\delta+1,2,1\}\pmod*{m}.
    \]
    
    \item If $\delta\equiv1\pmod*{m}$ and $\beta\equiv0\pmod*{m}$, then 
    \[
        L\equiv\{a,a+1,1,2\}\pmod*{m}
    \]
    and 
    \[
        R\equiv\{a+1,a,2,1\}\pmod*{m}.
    \]
    
    \item If $\delta\equiv-1\pmod*{m}$ and $\beta\equiv a\pmod*{m}$, then 
    \[
        L\equiv\{a,1,a-1,0\}\pmod*{m}
    \]
    and 
    \[
        R\equiv\{a-1,0,a,1\}\pmod*{m}.
    \]
    \end{enumerate}
    In each case, we have $L\equiv R\pmod*{m}$, as desired.
\end{proof}

\section{Explicit numerical semigroups and results}\label{sec:special-cases}
In the previous section, we focused on numerical semigroups for which had an Ap\'ery set whose nonzero elements formed an arithmetic sequence, and we derived explicit conditions to determine when such a numerical semigroup is equidistributed modulo $m$. In this section, we will explicitly describe the numerical semigroups that have such an Ap\'ery set.

As before, we suppose $S$ is a numerical semigroup with some nonzero $a\in S$ for which $\Ap(S;a)\setminus\{0\}$ is an arithmetic sequence. In other words, 
\[\Ap(S;a)=\{0\}\cup\{\beta+i\delta:i\in\ints{1,a-1}\}\]
for some $\beta,\delta\in\Z$. Recall that (by Lemma~\ref{lem:arith-props}) we have $\beta\ge0$, $\delta>0$, $a\mid \beta$, and $\gcd(a,\delta)=1$.

\subsection{Numerical semigroups for which the nonzero terms of an Ap\'ery set form an arithmetic sequence}
We begin by describing two well-known families of numerical semigroups. The first family consists of numerical semigroups generated by at most two integers.

\begin{prop}[{\cite[Section 3.I]{Selmer1977}}]\label{prop:emb-dim-2-apery-set}
    Suppose $S=\langle a,b\rangle$ with $\gcd(a,b)=1$. Then $\Ap(S;a)=\{0\}\cup\{ib : i\in\ints{1,a-1}\}$ and $\genus(S)=(a-1)(b-1)/2$. If $a=1$ or $b=1$, then $\embdim(S)=1$. Otherwise, $\embdim(S)=2$.
\end{prop}

The second family consists of numerical semigroups generated by \emph{generalized arithmetic sequences}, which are sequences of the form $\{a,ha+d,ha+2d,\dots,ha+kd\}$ for $a\ge2$ and $h,d,k\ge1$. In general, it is enough to consider $k\in\ints{1,a-1}$. For our work, we are interested in the case of maximal embedding dimension, which occurs when $k=a-1$.

\begin{prop}[{\cite[Section 3.III,IV]{Selmer1977}}]\label{prop:gen-arith-apery-set}
    Suppose $S=\langle\{a\}\cup\{ha+id : i\in\ints{1,a-1}\}\rangle$    with $a,h,d\in\N$ and $\gcd(a,d)=1$. Then $\embdim(S)=a$ and $\Ap(S;a)=\{0\}\cup\{ha+id : i\in\ints{1,a-1}\}$. Furthermore, $\genus(S)=(a-1)(2h+d-1)/2$.
\end{prop}

\begin{remark}\label{rmk:degenerate-case}
    When numerical semigroups generated by generalized arithmetic sequences appear in the literature, $h$ is taken to be a positive integer. (See, e.g., \cite{Lewin1975}, \cite{Selmer1977}, \cite{Ritter1998}, and \cite{Matthews2004}.) If we allow $h=0$, then 
    \[S=\langle\{a\}\cup\{id:i\in\ints{1,a-1}\}\rangle=\langle a,d\rangle,\]
    a numerical semigroup that is generated by at most two integers. We can think of the numerical semigroups with $h=0$ as a degenerate case of the family of numerical semigroups that are generated by generalized arithmetic sequences.
    
    For the remainder of this paper, we will work with the family of numerical semigroups of the form 
    \[S=\langle \{a\}\cup\{ha+id : i\in\ints{1,a-1}\}\rangle\] where $h\ge0$, $a,d\ge1$, and $\gcd(a,d)=1$. This family contains the two families described above (in Proposition~\ref{prop:emb-dim-2-apery-set} and Proposition~\ref{prop:gen-arith-apery-set}) as subfamilies.
\end{remark}
We now have a lemma which says Proposition~\ref{prop:emb-dim-2-apery-set} and Proposition~\ref{prop:gen-arith-apery-set} are ``if and only if'' statements. If we know an Ap\'ery set, then we know the semigroup.

\begin{lem}\label{lem:apery-iff}
    Let $S$ and $T$ be numerical semigroups with some nonzero $a\in S\cap T$. If $\Ap(S;a)=\Ap(T;a)$, then $S=T$.
\end{lem}
\begin{proof}
    Each element of $\Ap(S;a)$ is the minimal element of $S$ in its congruence class modulo $a$. Thus, 
    \[\HH(S)=\{w-ka : w\in\Ap(S;a),\, k\in\N,\, w-ka>0\}.\]
    Since $\Ap(S;a)=\Ap(T;a)$, $\HH(S)=\HH(T)$, and thus $S=T$.
\end{proof}

We are now able to explicitly describe the numerical semigroups where the nonzero terms of an Ap\'ery set form an arithmetic sequence.

\begin{prop}\label{prop:conds-for-two-families}
    Suppose $S$ is a numerical semigroup with some nonzero $a\in S$, and suppose $\Ap(S;a)\setminus\{0\}$ is an arithmetic sequence consisting of at least one term. Then $a\ge2$ and
    \[S=\langle\{a\}\cup\{ha+id : i\in\ints{1,a-1}\}\rangle\]
    for some $h\ge0$, $d\ge1$, and $\gcd(a,d)=1$. 
    
    If $h=0$ and $d=1$, then $\embdim(S)=1$. If $h=0$ and $d>1$, then $\embdim(S)=2$. If $h\ge1$, then $\embdim(S)=a$.
\end{prop}
\begin{proof}
     Since $\#(\Ap(S;a))=a$, if $\Ap(S;a)\setminus\{0\}$ contains at least one term, then $a\ge2$. We consider the cases of $a=2$ and $a\ge3$ separately.
     
     If $a=2$, then $S=\langle 2,b\rangle$ for some odd positive integer $b$. If $b=1$, then $S=\langle 2,1\rangle=\langle 1\rangle=\N_0$, in which case we have $\embdim(S)=1$, $h=0$, and $d=1$. If $b>1$, then by Proposition~\ref{prop:emb-dim-2-apery-set}, $\embdim(S)=2$ and we can write
     \[S=\langle 2,b\rangle = \langle\{2\}\cup\{0a+ib : i\in\ints{1,1}\}\rangle.\]
     Hence, we have $h=0$, $d=b>1$, and $\gcd(a,d)=\gcd(2,b)=1$.
     
     If $a\ge3$, and $\Ap(S;a)\setminus\{0\}$ is an arithmetic sequence, then 
     \[\Ap(S;a)\setminus\{0\}=\{\beta+i\delta : i\in\ints{1,a-1}\}\]
     for $\beta\ge0$, $\delta>0$, $a\mid \beta$, and $\gcd(a,\delta)=1$. Let $\beta= la$ for some $l\in\N_0$. We have two cases based on $l$. 
     
     If $l=0$, then
     \[\Ap(S;a)\setminus\{0\}=\{i\delta : i\in\ints{1,a-1}\}.\] 
     By Proposition~\ref{prop:emb-dim-2-apery-set} and Lemma~\ref{lem:apery-iff}, $S=\langle a,\delta\rangle$ with $\gcd(a,\delta)=1$. If $\delta=1$, then $S=\langle a,1\rangle=\langle 1\rangle=\N_0$, in which case we have $\embdim(S)=1$, $h=0$, and $d=\delta=1$. If $\delta>1$, then by Proposition~\ref{prop:emb-dim-2-apery-set}, $\embdim(S)=2$ and we can write
     \[S=\langle a,\delta\rangle
     =\langle\{a\}\cup\{0a+i\delta:i\in\ints{1,a-1}\}\rangle.
     \]
    Hence, we have $h=0$, $d=\delta>1$, and $\gcd(a,d)=\gcd(a,\delta)=1$.
     
     If $l>0$, then \[\Ap(S;a)\setminus\{0\}=\{la + i\delta : i\in\ints{1,a-1}\}.\]
     By Proposition~\ref{prop:gen-arith-apery-set} and Lemma~\ref{lem:apery-iff},
     \[S=\langle\{a\}\cup\{ha+id:i\in\ints{1,a-1}\}\rangle\]
     for $h=l$ and $d=\delta$. We have $h>0$, $d>0$, $\gcd(a,d)=1$, and $\embdim(S)=a$.
\end{proof}

\subsection{Explicit criteria for this family of numerical semigroups}
Now that we can explicitly describe the numerical semigroups $S$ for which the nonzero terms of an Ap\'ery set of $S$ form an arithmetic sequence, we can use the description to determine when the set of gaps of $S$ is equidistributed modulo $m$.
\begin{cor}\label{cor:gen-arith-seq-ED}
    Suppose $S$ is a numerical semigroup with some nonzero $a\in S$, and suppose $\Ap(S;a)\setminus\{0\}$ is an arithmetic sequence consisting of at least one term. Then $S=\langle \{a\}\cup\{ha+id:i\in\ints{1,a-1}\}\rangle$ for $a\ge2$, $h\ge0$, $d\ge1$, and $\gcd(a,d)=1$. We have that $\HH(S)$ is equidistributed modulo $m$ if and only if $\gcd(ad,m)=1$ and at least one of the following occurs: 
    \begin{enumerate}
        \item $a\equiv1\pmod*{m}$; or
        \item $a\equiv2\pmod*{m}$ and $2h+d\equiv1\pmod*{m}$; or
        \item $d\equiv1\pmod*{m}$ and $h\equiv0\pmod*{m}$; or
        \item $d\equiv-1\pmod*{m}$ and $h\equiv1\pmod*{m}$.
    \end{enumerate}
\end{cor}
\begin{proof}
     If $\Ap(S;a)\setminus\{0\}=\{\beta+i\delta:i\in\ints{1,a-1}\}$ is an arithmetic sequence consisting of at least one term, then by Proposition \ref{prop:conds-for-two-families}, we must have $S=\langle\{a\}\cup\{ha+id:i\in\ints{1,a-1}\}\rangle$ for $a\ge2$, $h\ge0$, $d\ge1$, and $\gcd(a,d)=1$. (As before, $\beta=ha$ and $\delta=d$.)
     
     By Theorem \ref{thm:main-result}, $\HH(S)$ is equidistributed modulo $m$ if and only if $\gcd(a\delta,m)=1$ and any of the four cases from Theorem \ref{thm:main-result} occur. For $n=1,2,3,4$, we will assume Case $n$ of Theorem~\ref{thm:main-result} and show that it is equivalent to Case $n$ of this corollary.
     
     $n=1$: Both cases have $a\equiv1\pmod*{m}$, so we are done.
     
     $n=2$: Suppose $a\equiv2\pmod*{m}$ and $\beta+\delta\equiv1\pmod{m}$. Since $\beta=ha\equiv2h\pmod*{m}$ and $\delta=d$, we equivalently have $a\equiv2\pmod*{m}$ and $2h+d\equiv1\pmod*{m}$. These two cases are equivalent.
     
     $n=3$: Suppose $\delta\equiv1\pmod*{m}$ and $\beta\equiv0\pmod*{m}$. Since $\beta=ha$ and $\delta=d$, we equivalently have $d\equiv1\pmod*{m}$ and $ha\equiv0\pmod*{m}$. The latter congruence is equivalent to the congruence $h\equiv0\pmod*{m}$ because $\gcd(a,m)=1$. These two cases are equivalent.
     
     $n=4$: Suppose $\delta\equiv-1\pmod*{m}$ and $\beta\equiv a\pmod*{m}$. Since $\beta=ha$ and $\delta=d$, we equivalently have $d\equiv-1\pmod*{m}$ and $ha\equiv a\pmod*{m}$. The latter congruence is equivalent to the congruence $h\equiv1\pmod*{m}$ because $\gcd(a,m)=1$. These two cases are equivalent.
    \end{proof}
    
    \begin{remark}\label{rmk:remark-revisited}
    We can revisit Proposition~\ref{prop:S-contains-3}, where we determined conditions for $S$ a numerical semigroup of multiplicity 3. In particular, when $\embdim(S)=\mult(S)=3$, we have $S=\langle 3,b,c\rangle$, and $\HH(S)$ is equidistributed modulo $m$ if and only if $\gcd(3,m)=1$ and $m$ divides either $\gcd(b-1,c-2)$ or $\gcd(b-2,c-1)$. In Remark~\ref{rmk:mult-emb-dim-3}, we note that the numerical semigroup $S=\langle 3,b,c\rangle$, with $b<c$, can be written as $S=\langle 3,3h+d,3h+2d\rangle$ for $h=(2b-c)/3$ and $d=c-b$. By Corollary~\ref{cor:gen-arith-seq-ED}, $\HH(S)$ is equidistributed modulo $m$ if and only if $\gcd(3d,m)=1$ and one of four conditions holds. As expected, working through the four conditions, we recover the same results that we obtained in Proposition~\ref{prop:S-contains-3}.
    \end{remark}

We can specialize to the case where $\embdim(S)=2$.

\begin{cor}\label{cor:emb-dim-2-ED}
    Suppose $\embdim(S)=2$, so $S=\langle a,b\rangle$ for $a,b>1$ with $\gcd(a,b)=1$. Then $\HH(S)$ is equidistributed modulo $m$ if and only if $\gcd(ab,m)=1$ and at least one of the following holds:
    \begin{enumerate}
        \item $a\equiv1\pmod*{m}$; or
        \item $b\equiv1\pmod*{m}$.
    \end{enumerate}
\end{cor}
\begin{proof}
    If $\embdim(S)=2$, then by Proposition~\ref{prop:emb-dim-2-apery-set} we have $S=\langle a,b\rangle$ for $a,b>1$ with $\gcd(a,b)=1$. Thus, $S=\langle \{a\}\cup\{ha+id:i\in\ints{1,a-1}\}\rangle$ for $h=0$ and $d=b$. Plugging these into the four cases of Corollary~\ref{cor:gen-arith-seq-ED}, we find that either $a\equiv1\pmod*{m}$ or $b\equiv1\pmod*{m}$. Hence $\HH(S)$ is equidistributed modulo $m$ if and only if $\gcd(ab,m)=1$ and at least one of $a$ and $b$ is 1 modulo $m$.
\end{proof}

Put another way, for $S=\langle a,b\rangle$, $\HH(S)$ is equidistributed modulo $m$ if and only if $m$ is a divisor of $a-1$ with $\gcd(b,m)=1$, or if $m$ is a divisor of $b-1$ with $\gcd(a,m)=1$. (Note that we can allow $a=1$ or $b=1$ here as well.)

\begin{example}\label{ex:ED-example-revisited}
    Let $S=\langle 5,7\rangle$. Then $\HH(S)=\{1,2,3,4,6,8,9,11,13,16,18,23\}$. In Example~\ref{ex:ED-example}, we determined that this set is equidistributed modulo $m$ for $m\in\{1,2,3,4,6\}$. We will revisit this example now that we have Corollary~\ref{cor:emb-dim-2-ED}. The divisors of $5-1$ are each relatively prime to 7. The divisors of $7-1$ are each relatively prime to 5. Thus, $\HH(S)$ is equidistributed modulo any divisor of 4 or 6, which is to say $\HH(S)$ is equidistributed modulo 1, 2, 3, 4, and 6.
\end{example}

For a numerical semigroup $S$ of maximal embedding dimension generated by a (purely) arithmetic sequence, we have $S=\langle a, a+d, a+2d, \dots, a+(a-1)d\rangle$ for $a,d\ge1$. We can obtain results for such a numerical semigroup by plugging $h=1$ into Corollary \ref{cor:gen-arith-seq-ED}. Note that the third case can now only occur if $m=1$, and the second case is now a special case of the fourth case. We therefore have two cases: $a\equiv1\pmod*{m}$ or $d\equiv-1\pmod*{m}$.

\begin{cor}\label{cor:pure-arith-ED}
    Suppose $S$ is a numerical semigroup of maximal embedding dimension generated by an arithmetic sequence. Then $S=\langle \{a+id:i\in\ints{0,a-1}\}\rangle$ for $a\ge1$, $d\in\N$, and $\gcd(a,d)=1$. We have that $\HH(S)$ is equidistributed modulo $m$ if and only if $\gcd(ad,m)=1$ and at least one of the following occurs: 
    \begin{enumerate} 
        \item $a\equiv1\pmod*{m}$; or 
        \item $d\equiv-1\pmod*{m}$.
    \end{enumerate}
\end{cor}

\subsection{Observations}
Finally, we make a few observations.

For $\embdim(S)=2$, $\genus(S)=(a-1)(b-1)/2$. By Corollary~\ref{cor:emb-dim-2-ED}, $\HH(S)$ is equidistributed modulo $m$ if and only if $\gcd(ab,m)=1$ and either $m\mid(a-1)$ or $m\mid(b-1)$.

In the purely arithmetic case, $\genus(S)=(a-1)(d+1)/2$. By Corollary~\ref{cor:pure-arith-ED}, $\HH(S)$ is equidistributed modulo $m$ if and only if $\gcd(ad,m)=1$ and either $m\mid(a-1)$ or $m\mid(d+1)$.

These two cases have a very similar feel. The generalized arithmetic case is similar, though we don't quite get an ``if and only if'' result.

In the generalized arithmetic case, $\genus(S)=(a-1)(2h+d-1)/2$. By Corollary~\ref{cor:gen-arith-seq-ED}, if $\HH(S)$ is equidistributed modulo $m$, then $\gcd(ad,m)=1$ and either $m\mid(a-1)$ or $m\mid(2h+d-1)$, the latter congruence implied by, but not equivalent to, cases 2, 3, 4 of Corollary~\ref{cor:gen-arith-seq-ED}.

\section*{Acknowledgments}
The author wishes to thank the anonymous referees for their detailed reports. Their suggestions were incredibly helpful. The author also wishes to thank Pieter Moree for suggesting the phrase `equidistribution modulo $m$'.

\nocite{*}
\bibliographystyle{abbrvnat}
\bibliography{refs}
\end{document}